\documentclass[12pt]{article}

\usepackage[T2A]{fontenc}
\usepackage[cp1251]{inputenc}
\usepackage[russian]{babel}

\usepackage{amssymb}
\usepackage{amsmath}
\usepackage{euscript}

\begin{document}

\centerline{\large\bf Inertial Manifolds and Limit Cycles}
\centerline{\large\bf of Dynamical Systems in ${\mathbb R}^{n}$}

\bigskip

\centerline{\bf L.A.~Kondratieva$^{1}$, A.V.~Romanov$^{2}$}

\bigskip
\noindent $^{1}$ Moscow Aviation Institute (National Research
University), Russia

\noindent $^{2}$ National Research University Higher School of
Economics, Moscow, Russia,

\textit{Email adress}: av.romanov@hse.ru

\bigskip

\textbf{Abstract}. We show that the presence of a two-dimensional
inertial manifold for an ordinary differential equation
in~${\mathbb R}^{n}$ permits reducing the problem of determining
asymptotically orbitally stable limit cycles to the
Poincare--Bendixson theory. In the case $n=3$ we implement such a
scenario for a model of a satellite rotation around a celestial
body of small mass and for a biochemical model.

\medskip



\bigskip

\noindent \textbf{1. Introduction}

\noindent We consider ordinary differential equations
$$
\dot{x}=-Ax+F(x),\quad x\in {\mathbb R}^{n} ,\quad n\ge 3,
\eqno (1.1)
$$
where $A$ is a symmetric $n\times n$ matrix with eigenvalues $0\le
\lambda_{1} \le \lambda_{2} \le \dots\le \lambda_{n}$ and the
function $F$ belongs to
$C^{1+\alpha}(\mathbb{R}^{n},\mathbb{R}^{n})$ for some $\alpha\in
(0,1)$. We let $F'(x)$ denote the Jacobi matrix of the mapping $F$
at a point~$x$, and $\left\|\cdot\right\|$ and $\left\|
\cdot\right\|_{2}$ denote the Euclidean norm in ${\mathbb R}^{n} $
and the Euclidean norm of matrices, respectively. If one of the
two conditions
$$
\left\|F(x)-F(y)\right\| \le K\left\| x-y\right\| ,  \; \;
\left\|F'(x)\right\|_{2} \le K,\quad x,y\in {\mathbb R}^{n},
\eqno(1.2)
$$
that are equivalent in this situation is satisfied, then equation
(1.1) generates a $C^{1}$-smooth phase flow
$\{\Phi_{t\in\mathbb{R}}\}$ in ${\mathbb R}^{n}$. Everywhere below
we identify linear operators on ${\mathbb R}^{n}$ with their
matrices.
  Let $f=-A+F$ be a vector field of
(1.1), then we call $x_{s}\in \mathbb{R}^{n}$ a \textit{singular
point} if $f(x_{s})=0$. By a \textit{cycle} we mean a closed
trajectory. A \textit{stable limit cycle} is a cycle that is
asymptotically orbitally stable as $t\to +\infty$.

The theory of inertial (that is, invariant and globally
exponentially attracting) manifolds was developed in the 1980s as
a tool for studying the final (at large times) dynamics of
semilinear parabolic equations with a vector field structure of
the form~(1.1) in an \textit{infinite-dimensional} Hilbert
space~$X$ (see [1, Ch.~8], [2] and the references therein). In
this case, as usual, it is assumed that~$A$ is an unbounded
self-adjoint positive linear operator in~$X$ with a compact
resolvent. In such a situation, the presence of an $m$-dimensional
inertial manifold (IM) permits describing the final dynamics of an
infinite-dimensional evolutionary system by an ordinary
differential equation (ODE) in ${\mathbb R}^{m}$. Here we
demonstrate the usefulness of inertial manifolds in the
finite-dimensional case $X={\mathbb R}^{n}$. Namely, the existence
of a two-dimensional IM $(m=2)$ allows one to reduce studying the
final dynamics of equation (1.1) to solving the corresponding
problem in~${\mathbb R}^{2}$ and, in several cases, to prove the
presence and to discover the localization of a stable limit cycle
without using the bifurcation technique or some rather complicated
topological constructions. We stress that, in contrast to the
bifurcation theory, our approach proves the existence of stable
self-sustained oscillations of a ``large amplitude''.

\bigskip

\noindent \textbf{2. Inertial manifolds}

\noindent A set $\Lambda \subseteq {\mathbb R}^{n}$ is said to be
invariant if $\Phi_{t} \Lambda =\Lambda $, $t>0$. Let $P_{m}$ and
$Q_{m}$ be orthogonal projection operators in ${\mathbb R}^{n}$ on
the subspaces $X_{m}$ and $X_{n-m}$ corresponding to the
eigenvalues $\lambda_{1} ,\dots,\lambda_{m} $ and $\lambda_{m+1}
,\dots,\lambda_{n}$, $\lambda_{m}<\lambda_{m+1}$, of the
matrix~$A$.

Invariant manifold of the form
$$
H_{m} =\{ x\in {\mathbb R}^{n} :\; x=u+h(u),\;\; u\in X_{m}\}
\eqno (2.1)
$$
with the function $h\in \mathrm{Lip}\,(X_{m},X_{n-m})\bigcap
C^{1}(X_{m},X_{n-m})$ we call inertial, if for each trajectory
$x(t)$, there exists a trajectory $\overline{x}(t)\subset H_{m}$
such that
$$
\left\| x(0)-\overline{x}(0)\right\| \le M_{1}
\left\|Q_{m}x(0)-h(P_{m}x(0))\right\|, \eqno (2.2)
$$
$$
\left\| x(t) -\, \overline{x}(t)\right\| \; \le \; M_{2}e^{-\gamma
t} \left\| x(0) - \,\overline{x}(0)\right\| \eqno(2.3)
$$
for $t>0$, where $M_{1},M_{2}, \gamma>0$. If a set $E\subset
\mathbb{R}^{n}$ is bounded, then the Lipschitzian function
$h:X_{m}\rightarrow X_{n-m}$ is bounded on the bounded set
$P_{m}E$ and for everyone $x(0)\in E$ we have
$\|Q_{m}x(0)-h(P_{m}x(0))\|\leq M$ with $M=M(E)$. It follows from
(2.3) that $\|x(t)-\,\overline{x}(t)\| \le M_{1}M_{2}Me^{-\gamma
t}$ for $x(0)\in E,\,t>0$, which means $H_{m}$ exponentially and
uniformly attracts $E$. Let $\Lambda\subset \mathbb{R}^{n}$ be a
compact invariant set and $y\in \Lambda$. If $x(0)=\Phi_{-t}y$,
then $x(0)\in\Lambda,\;x(t)=y$, and
$$\|x(t) -\, \overline{x}(t)\|=\|y-\,\overline{x}(t)\|\le \; \overline{M}(\Lambda)e^{-\gamma t}.$$
Since $t>0$ is arbitrary, $\overline{x}(t)\in H_{m}$ and the set
$H_{m}$ is closed, then $y\in H_{m}$ and $\Lambda \subset H_{m}$.
In this way, the inertial manifold contains all compact invariant
sets (including the singular points and cycles) of the dynamical
system.

It is well known [3, 4] that if the \textit{exact spectral gap
condition}
$$
\lambda_{m+1} \, - \, \lambda_{m} \,  >\,  2K \eqno (2.4)
$$
is satisfied, then there is such a manifold with $h\in
\mathrm{Lip}\,(X_{m},X_{n-m})$ and the factor~$2$ on the
right-hand side of (2.3) cannot be decreased in general. Later, it
was shown~[2], that condition (2.4) also provides the existence of
a $C^{1}$-smooth inertial manifold. Estimate~(2.2) means that
$\left\| x(0) - \overline{x}(0)\right\|$ is small if the initial
point $x(0)$ is close to $H_{m}$. Estimate (2.3) reflects the
exponential tracking of the initial trajectory $x(t)$ by the
trajectory $\overline{x}(t)\subset H_{m}$.

By the reduction principle [4, Lemma~1], the compact invariant
sets $\Lambda$ of equation (1.1) and $P_{m}\Lambda$ of the ODE
$$
\dot{u}=-Au+P_{m}F(u+h(u)),\quad \; u=P_{m}x, \eqno (2.5)
$$
in $X_{m} \simeq {\mathbb R}^{m}$ are simultaneously
asymptotically stable or unstable. The dynamical system generated
by~(2.5) is topologically conjugate to the restriction of the
original dynamical system~(1.1) to~$H_{m}$. This means that the
final (for $t \rightarrow + \infty $) regimes of the original
equation in $\mathbb{R}^{n} $ are fully described by some ODE in
space of smaller dimension, which in many cases simplifies their
research. Essentially, we highlight the $m <n$ ``defining''
degrees of freedom of a $n$-dimensional dynamical system. In
addition, if $t$ is sufficiently large then every solution $x(t)$
of equation (1.1) is completely determined by its projection
$u(t)=P_{m}x(t)$ onto the subspace $X_{m}$ and is reconstructed by
the formula $x(t)=\psi (u(t))$ with $\psi (u)=u+h(u)$.

Splitting the right-hand side of equation (1.1) into linear and
nonlinear components, of course, is not unique. Right choice
matrix $A$ in (1.1) can help to satisfy the condition (2.4). On
the other hand, condition (2.4) can sometimes be ensured by using
a nondegenerate linear change of variables; the topology of the
phase portrait of the dynamical system does not change in this
case. Such a method is used below in Section~4 to study a
mathematical model of cell processes.

\noindent \textbf{Remark~2.1.} The existence of a two-dimensional
inertial manifold allows one to assert that the union of all
singular points and cycles (if any) has the form of a Lipschitz
graph over a certain plane $X_{2}\subset{\mathbb R}^{n}$.

It should be noted that, under condition~(2.4), the inertial
manifold $H_{m}$ does not inherit the smoothness of the
nonlinearity~$F$; for example, the condition that~$F$ is real
analytic in~$\mathbb{R}^{n}$ does not even imply that $H_{m}\in
C^{2}$.

\noindent \textbf{Definition~2.1.} A domain $D\subset {\mathbb
R}^{n}$ is strictly positive invariant if
$\Phi_{t}\overline{D}\subseteq D$, $t>0$.

In particular, this means that the boundary $\partial D$ does not
contain singular points.

\noindent \textbf{Remark 2.2.} Even under a weaker condition
$\Phi_{t} D\subseteq D$, $t>0$, the continuity of the mapping
$x\rightarrow \Phi_{t}x$ for $x\in {\mathbb R}^{n}$ guarantees the
inclusion $\Phi_{t}\overline{D}\subseteq \overline{D}$, $t>0$, for
the closure~$\overline{D}$.

The strict positive invariance of $D$ is ensured if the vector
field $f(x)=-Ax+F(x)$ of equation (1.1) on the boundary $\partial
D$ is directed inside the interior of~$D$. If the domain $D\subset
{\mathbb R}^{n}$ is strictly positive invariant, then the domain
$P_{m}D\subset X_{m}$ has the same property with respect to the
ODE~(2.5).

\noindent \textbf{Remark 2.3.} The closure of the union of all
cycles contained in the strictly positive invariant domain~$D$
does not contain points of~$\partial D$.

This is a consequence of the continuity of the phase
flow~$\{{\Phi_{t}}\}$ with respect to $x\in \mathbb{R}^{n}$.

Consider the quadratic form $V(x)=\left\|
Qx\right\|^{2}-\left\|Px\right\|^{2}$ with an arbitrary orthogonal
projection operator~$P$ in~$\mathbb{R}^{n}$ and $Q=\mathrm{Id}-P$.
Assume that, for some $\lambda ,\varepsilon>0$, any two solutions
$x(t)$ and $y(t)$ of (1.1) satisfy the following relation holds
with $t>0$:
$$
\frac{d}{dt} V(x(t)-y(t))+2\lambda V(x(t)-y(t))\le -\varepsilon
\left\|x(t)-y(t)\right\|^{2}. \eqno (2.6)
$$
This condition is known in the theory of inertial manifolds as the
\textit{strong cone condition}.

\noindent \textbf{Remark 2.4} (see [2, Lemma~2.21; 4, Lemma 4]).
Condition (2.4) implies (2.6) with $P=P_{m},\;\lambda
=(\lambda_{m+1} +\lambda_{m} )/2$ and
$\varepsilon=(\lambda_{m+1}-\lambda_{m})/2-K$.

Recall the well-known (see [5]) estimate
$\mathcal{T}\geq2\pi/K_{1}$ of the periods $\mathcal{T}>0$ of
periodic solutions~(1.1), where $K_{1}=\lambda_{n}+K$ is the
Lipschitz constant of the vector field $f=-A+F$. For $\tau
=\pi/K_{1}$, we set
$U_{\tau}(x)=x-\Phi_{\tau}x,\;x\in\mathbb{R}^{n}$. The zeros of
the vector field $U_{\tau}$ are precisely the singular points of
equation (1.1). A  point $x_{s}$  is said to be asymptotically
unstable if the spectrum $\sigma(f'(x_{s}))$ contains an
eigenvalue with $\mathrm{Re}\,\lambda>0$. In this case, $\sigma
(U'_{\tau}(x_{s}))=\{1\}- \mathrm{exp}(\tau\sigma(f'(x_{s})))$.

\noindent \textbf{Theorem~2.1.} \textit{Assume that the following
conditions are satisfied for equation} (1.1):

\noindent (i) \textit{there exists bounded convex strictly
positive invariant domain $D\subset {\mathbb R}^{n}$ containing a
unique singular point~$x_{s}$, this point  is asymptotically
unstable and satisfies} $\mathrm{det}f ^{{'}}(x_{s})\neq 0$;

\noindent (ii) \textit{the function $F$ is real analytic in~$D$};

\noindent (iii) $\lambda_{3}-\lambda_{2}>2K$.

\textit{Then at least one stable limit cycle is localized in the domain}~$D$.

\noindent \textit{Proof}. We use condition~(iii) to reduce the
final dynamics of~(1.1) to the two-dimensional inertial manifold
$H_{2}\ni x_{s}$. By Remark~2.4, the estimate~(iii) implies
relation~(2.6) for the quadratic form~$V$ with
$P=P_{2},\;\lambda=(\lambda_{3} +\lambda_{2} )/2$ and
$\varepsilon=(\lambda_{3}-\lambda_{2})/2-K$. Assume that
$\mathrm{Re}\,\kappa_{1} \ge \mathrm{Re}\,\kappa_{2} \ge \dots\ge
\mathrm{Re}\,\kappa_{n}$ for $\kappa_{i}\in \sigma(f'(x_{s}))$. If
we consider the matrix $f'(x_{s})$ as a perturbation of the
matrix~$-A$, then condition~(iii) implies the inequality
$\mathrm{Re}\,\kappa_{3}<-\lambda<0$. It follows from
condition~(i) that the vector field $U_{\tau}$ with
$\tau=\pi/K_{1}$ has a unique zero~$x_{s}$ in~$\overline D$.

Since the domain $D$ is convex and $\Phi_{\tau}\overline {D}
\subset D$, then according to [6, Theorem 21.5] the vector field $
U_{\tau}$ is not is degenerate (0 does not belong to $\sigma (U'_
{\tau})$) on $\partial D$ and the rotation of $U_{\tau}$ on
$\partial D$ is equal to 1. By the hypothesis (i) of the theorem
the vector field $U_{\tau}$ is not degenerate at the point
$x_{s}$, therefore from [6, Theorem 20.6] and [6, Theorem 21.6] we
successively find that $\mathrm{ind}\, x_{s}=1$ and $\mathrm
{ind}\, x_{s} = (- 1)^{\beta}$, where $\mathrm {ind}$ is the
Poincare index and $\beta$ is an even sum multiplicities of the
real $\lambda> 1 $ in $\sigma (\Phi'_{\tau}(x_{s}))$. At the same
time, $\beta$ is the sum multiplicities of positive $\kappa \in
\sigma(f'(x_{s}))$. So, since $\mathrm {Re}\, \kappa_{3} <0 $ and
$\mathrm {Re}\, \kappa_ {1}>0$, then $\mathrm {Re} \,
\kappa_{2}>0$.

Thus, taking ~(i), (ii), and Remark~2.3 into account, we see that
the assumptions in~[7, Corollary~6.1] are satisfied, and hence the
domain~$D$ contains at most finitely many cycles. One can see that
the point $P_{2}x_{s}$ is an unstable focus or an unstable knot of
equation~(2.5) in the plane $X_{2} \subset {\mathbb R}^{n}$. By
the Poincare--Bendixson theory [8, Sect.~2.8], this equation has
finitely many embedded cycles in the strictly positive invariant
domain $P_{2}D\subset X_{2}$ and at least one of them, $\Gamma$,
is stable. Then $\psi \Gamma$ is a stable limit cycle of the
original equation~(1.1). $\Box$

Theorem~2.1 gives us a method for determining stable limit cycles
of ODEs in~${\mathbb R}^{n}$. In what follows we refer to this
method as to the ``spectral gap method''. In fact, notion similar
to that of inertial manifold has been used successfully by
R.A.~Smith (see [7, 9, 10] and the references therein) in his
studies of cycles of ODEs. This author worked with Lipschitz
invariant manifolds of the form (2.1), attracting (not necessarily
exponentially) all trajectories for $t\rightarrow+\infty$ and
containing all bounded invariant sets. He did not use the simple
and convenient condition~(2.4) but directly
considered\footnote{\,See, e.g., [10, Theorem~3].} the condition
of type (2.6) with an arbitrary quadratic form $V(x)$ of the
signature $(0,n-2,2)$. Formally, assumption~(2.6) is weaker
than~(2.4) and does not mean that the vector field of the equation
splits into linear and nonlinear parts. At the same time, the
spectral gap condition~(2.4) can be verified significantly
simpler.

On the other hand, the method proposed in~[3] guarantees the
existence of an inertial manifold of dimension $m<n$ for equations
of the form~(1.1) with an \textit{arbitrary} linear part~$-A$ if,
for some $\lambda>0$, the spectrum $\sigma(A)$ has~$m$ values
(with multiplicity taken into account) in the half-plane
$\mathrm{Re}\,z<\lambda$, the straight line
$\mathrm{Re}\,z=\lambda$ lies in the resolvent set $\rho(A)$, and
$\|(A-\lambda-i\omega)^{-1}\|_{2}<1/K$, $\omega\in \mathbb{R}$.
Such a technique was independently used to determine stable limit
cycles in~[10]. The author believes that the revival of this
approach is rather perspective.

It should be noted that the technique of this paper (as well as
papers [7, 9, 10]) only detect ODE cycles lying on invariant
$2\mathrm{D}$-manifolds of the Cartesian structure (2.1).

In the following two sections we illustrate the spectral gap
method with examples from two distinct areas of natural science.

\bigskip
\bigskip
\bigskip

\noindent \textbf{3. Satellite motion model}

\noindent The problems of the periodic dynamics of the satellites
of celestial bodies extensive literature is devoted (see, for
example, [11] and references therein). In particular, the dynamics
of a artificial satellite flying around a celestial body of small
mass was studied in [12]. We consider here this model as a
successful mathematical application of our method for detecting
stable limit cycles. Let $(r,\varphi )$ be the polar coordinates
in the plane of the motion $r=r(t),\;\varphi=\varphi(t)$ of a
flying vehicle. According to [12], the radial and transverse
control forces act on the satellite, depending on the positive
parameters $\mu_{1},\mu_{2},\mu_{3}$ and some smooth function
$g(\dot{\varphi})$. The goal is to determine the values
$\mu_{1},\mu_{2},\mu_{3}$ and the function $g$ so as to ensure the
existence of a stable periodic motion in coordinates
$(r,\dot{r},\dot{\varphi})$. We set $x_{1} =\dot{r}+\mu_{2} r$,
$x_{2} =r$, $x_{3} =\dot{\varphi}$. In these new coordinates, the
satellite dynamics can be described by the system of equations
(slightly different from the system in~[12])
$$
\quad \; \dot{x}_{1}=-\mu_{1} x_{1}+g(x_{3}),
$$
$$
\dot{x}_{2}=-\mu_{2} x_{2} +x_{1},
$$
$$
\dot{x}_{3}=-\mu_{3} x_{3} +x_{2} \eqno (3.1)
$$
with control parameters $\mu_{1},\mu_{2},\mu_{3}>0$ and the
``admissible'' nonlinear function $g\in C^{1+\alpha}({\mathbb
R})$. We define the class of admissible smooth functions $g$ in
(3.1) by conditions
$$
0<g(x_{3})\, <M,\quad \; -1\, \le \, g'(x_{3})\, < \, 0
\eqno (3.2)
$$
for $x_{3}\in {\mathbb{R}}$. The choice of such a class will allow
us to apply Theorem 2.1 under certain conditions on the parameters
$\mu_{1},\mu_{2},\mu_{3}$. A similar mathematical model was
studied in~[10, Sect.~7] from a different standpoint. System~(3.1)
takes the form~(1.1) if we set
$$
A=\left(\begin{array}{ccc} {\mu_{1} }
& {0} & {0} \\ {0} & {\mu_{2} } & {0} \\ {0} & {0} & {\mu_{3} }
\end{array}\right)\;
{\kern 1pt}, \;\;\; F(x)=\; \, \left(\begin{array}{c} {g(x_{3} )} \\ {x_{1} } \\ {x_{2} }
\end{array}\right).
$$
This decomposition of a vector field (3.1) is natural from the
point of view of condition (iii) of Theorem 2.1, so as the matrix
$A$ is symmetric, and the Lipschitz constant of nonlinearity $F$
easy to appreciate.

Due to the second condition in (3.2), system~(3.1) generates a
$C^{1}$ phase flow $\{\Phi_{t}\}$ in~${\mathbb R}^{3}$.

\noindent \textbf{Lemma~3.1.} \textit{The convex domain}
$$
D=\{x\in{\mathbb R}^{3} :\;  0<x_{1} <\; \frac{M}{\mu_{1} } ,\; \;
0<x_{2} <\; \frac{M}{\mu_{1} \mu_{2} } ,\; \; 0<x_{3} <\;
\frac{M}{\mu_{1} \mu_{2} \mu_{3}}\}
$$
\textit{is strictly positive invariant and contains a unique
singular point}.

\noindent \textit{Proof}. The search of the singular points of the
system (3.1) reduces to solving the scalar equation
$g(x_{3})=\mu_{1} \mu_{2} \mu_{3} x_{3}$. Since according to
conditions (3.2) we have $0<g<M$ and $g'<0$, then this equation
has a unique solution $x_{3} =\nu
>0$. So there exists a unique singular point in~${\mathbb
R}^{3}$:
$$
x_{s} =\left(\mu_{2} \mu_{3} \nu ,\mu_{3} \nu ,\nu\right)
=\left(\frac{g(\nu )}{\mu_{1} } ,{\kern 1pt} \,\frac{g(\nu
)}{\mu_{1} \mu_{2} } ,\, \frac{g(\nu )}{\mu_{1} \mu_{2} \mu_{3} }
\right).
$$
Note that $x_{s} \in D$.

We first show that $\Phi_{t} D\subseteq D$, and hence
$\Phi_{t}\overline{D}\subseteq \overline{D}$ for $t>0$. Consider
the solution $x(t)=(x_{1}(t),x_{2}(t),x_{3}(t)$ with $x(0)\in D$.
On the faces $x_{1} =0$ and $x_{1} = M/\mu_{1}$ of the
parallelepiped~$D$, we have $\dot{x}_{1} =g(x_{3} )>0$ and
$\dot{x}_{1} =-\mu x_{1} +g(x_{3} )<0$ respectively, so
$0<x_{1}(t)<M/\mu_{1}$ for $t>0$. On the faces $x_{2} = M/(\mu_{1}
\mu_{2} )$ and $x_{2} =0$, we have $\dot{x}_{2} <0$ and
$\dot{x}_{2} (t)=x_{1}(t)>0$ respectively, and hence,
$0<x_{2}(t)<M/(\mu_{1}\mu_{2})$ for $t>0$. On the faces $x_{3} =
M/(\mu_{1} \mu_{2} \mu_{3} )$ and $x_{3}=0$, we have $\dot{x}_{3}
<0$ and $\dot{x}_{3} (t)=x_{2} (t)>0$ respectively, so that
$0<x_{3}(t)<M/(\mu_{1}\mu_{2}\mu_{3})$ for $t>0$.

We write $\Pi =\{ x\in \partial D:\Phi_{t} x\in D,\;t>0\}$ and
$\Pi_{0} =\partial D\backslash \Pi$. We see that $\Pi_{0}
\subseteq l_{1} \bigcup l_{2} \bigcup \{ 0\} $, where $l_{1} =\{
x\in \partial D:x_{1} =0,\, x_{2}=0,\, x_{3}>0\} $ and $l_{2} =\{
x\in \partial D:x_{1} >0,\, x_{2} =0,\, x_{3}=0\}$. On $l_{1}$ and
$l_{2}$, we respectively have $\dot{x}_{1} >0$ and $\dot{x}_{2}
>0$, and hence $\Phi_{t} x\in D$, $t>0$, on $\Pi_{0} /\{ 0\} $.
Because $\Phi_{t} 0\ne 0$, we have $\Phi_{t} 0\in D$, $t>0$. Thus,
$\Pi_{0} =\phi$, $\Pi =\partial D$, and $\Phi_{t}
\overline{D}\subseteq D$ for $t>0$. $\Box$

Clearly,

$F'(x)=\left(\begin{array}{ccc} {0} & {0} & {g'(x_{3} )}
\\ {1} & {0} & {0} \\ {0} & {1} & {0} \end{array}\right),\quad \;
(F'(x))^{*} \cdot F'(x)={\rm diag}\, {\rm (1,}\, {\rm 1,}\,
(g'(x_{3} ))^{{\rm 2}} {\kern 1pt} {\rm )},$

\noindent and $\left\| F'(x)\right\|_{2}=1$ for all $x\in{\mathbb
R}^{3}$. Let $\lambda_{1}, \lambda_{2}, \lambda_{3}$ stand for the
parameters $\mu_{1}$, $\mu_{2}$, $\mu_{3}$ permutated by
nondecreasing order. We have $K=1$ and the spectral gap
condition~(2.4) becomes
$$
\lambda_{3} -\lambda_{2} >2. \eqno (3.3)
$$
We linearize the vector field of the system (3.1) at the singular
point $x_{s}$. Note that the Routh--Hurwitz criterion gives the
condition of asymptotic instability of $x_{s}$ by the inequality
$$
-g'(\nu )\, +\, \, \lambda_{1} \lambda_{2} \lambda_{3} \, \, > \,
\, \left(\lambda_{1} +\lambda_{2} +\lambda_{3}\right)
\left(\lambda_{1} \lambda_{2} +\lambda_{1} \lambda_{3}+\lambda_{2} \lambda_{3} \right).
\eqno (3.4)
$$
In addition,
$\mathrm{det}(F^{{'}}(x_{s})-A)=g'(\nu)-\lambda_{1}\lambda_{2}\lambda_{3}\neq
0$. Estimates (3.3), (3.4) determine a nonempty open set~$\Omega$
in the positive octant ${\mathbb R}_{+}^{3}$ of the parameters
$\lambda_{1}$, $\lambda_{2}$, $\lambda_{3}$. In particular, the
domain~$\Omega$ contains points of the form ($\delta,\; \delta,\;
2+2\delta$) for all sufficiently small $\delta
>0$. If the function $g$ in~(3.2) is real analytic for
$0<x_{3}<M/(\mu_{1}\mu_{2}\mu_{3})$, then by Theorem~2.1,
system~(3.1) with $(\lambda_{1} ,\lambda_{2} ,\lambda_{3}
)\in\Omega$ has a stable limit cycle $\Gamma \subset D$.

As an admissible nonlinear function in (3.1) we can, for example,
take
$$
g(x_{3} )={\rm arccot}(x_{3} -\nu ),\quad \nu=\frac{\pi }{2\mu_{1}
\mu_{2} \mu_{3}}.
$$
This function satisfies conditions (3.2) with $g'(\nu )=-1$ and
$M=\pi$.

In similar constructions~[12], the real analyticity of the
function~$g$ in (3.1) is not required, but it is only necessary to
prove the existence of an orbitally stable periodic trajectory on
which at least one different trajectory is  ``winding'' as
$t\rightarrow +\infty$.

\bigskip

\noindent \textbf{4. A model of cell processes}

\noindent Another example illustrating the spectral gap method is
related to the complex dynamics in cell processes~[13]. Consider
the following the system of equations
$$
\;\dot{x}=-kx+R(z),
$$
$$
\dot{y}=x-G(y,z),
$$
$$
\;\;\,\dot{z}=-qz+G(y,z),
\eqno (4.1)
$$
where
\[
R(z)=\frac{1}{1+z^{4} }, \quad G(y,z)=\frac{Ty(1+y)(1+z)^{2}}{L+(1+y)^{2} (1+z)^{2} }
\]
and $k,q, T, L>0$ are constants. Here $x$, $y$, and $z$ are
dimensionless concentrations of the matters $S_{1}$, $S_{2}$, and
$S_{3}$, where $S_{1}$ is the initial product, $S_{2}$ is the
intermediate product, and $S_{3}$ is the final product; $k$ and
$q$ are constants of the rate of variation in $S_{1}$ and $S_{3}$.
We have
\[
R_{z} =-\frac{4z^{3}}{(1+z^{4})^{2} },\;\; G_{z} =\frac{2TLy(1+y)(1+z)}{(L+(1+y)^{2} (1+z)^{2} )^{2} } \, ,
\]
\[
G_{y} =\frac{2TLy(1+z)^{2} }{(L+(1+y)^{2} (1+z)^{2} )^{2} }
+\frac{T(1+z)^{2} }{L+(1+y)^{2} (1+z)^{2} }\,,
\]
$R_{z} (z)<0$ for $z>0$, and $G(y,z)<T$, $G_{y}(y,z)>0$,
$G_{z}(y,z)>0$ for $y,z>0$. Since the first derivatives of the
functions $R$ and $G$ are uniformly bounded in $z\in{\mathbb R}$
and $(y,z)\in {\mathbb R}^{2}$, we see that system (4.1) generates
a smooth flow $\{\Phi_{t}\}$ in ${\mathbb R}^{3}$. We fix the
values $T=10$ and $L=10^{6}$ that are physically meaningful from
the standpoint of the authors of~[13] and try to determine pairs
of free parameters $(k,q)\in \mathbb{R}_{+}^{2}$ for which this
system satisfies the conditions of Theorem~2.1 and hence admits a
stable periodic regime.

Everywhere below we restrict ourselfs to the simple case when
$kT>1$ and $k>q$. By $p(x,y,z)$ we denote points in ${\mathbb
R}^{3}$.

\noindent \textbf{4.1. Positive invariant domain and a singular
point}. We note that $G(+\infty ,0)=T$ and $G(0,z)=0$ for $z>0$.
Since $kT>1$, we can uniquely determine the value $y_{0}>0$ from
the relation $G(y_{0} ,0)=1/k$. In what follows we set
$x_{0}=1/k,\;z_{0}=T/q$.

\noindent \textbf{Lemma~4.1.} \textit{The convex domain
$D=\{p\in{\mathbb R}^{3}: 0<x<x_{0},\; 0<y<y_{0},\; 0<z<z_{0}\}$
is strictly positive invariant and contains a unique singular
point}.

\noindent \textit{Proof}. Equating the right-hand side of~(4.1) to
zero we obtain the relations $x=qz$ and $kqz=R(z)$ which are
satisfied for a unique pair of values $x_{s},\, z_{s} >0$. Another
scalar equation $\varphi (y)=0$ with $\varphi (y)=qz_{s}-G(y,z_{s}
)$, $\varphi'<0$, has a unique solution $y_{s}>0$. So system~(4.1)
has a unique singular point $p_{s} =(x_{s},\,y_{s},\,z_{s})$ in
${\mathbb R}_{+}^{3}$. Since the function~$R$ decreases in $z>0$,
it follows that $z_{s} =(kq)^{-1} R(z_{s} )<(kq)^{-1}<z_{0}$ and
$x_{s} =k^{-1}R(z_{s} )<x_{0}$. Taking into account that $G$ is an
increasing function with respect to each variable $y>0$ and $z>0$,
from the relation $x_{s} =G(y_{s},z_{s} )$ we derive that
$x_{s}=G(y_{s},z_{s}) <x_{0}=G(y_{0},0)$, and hence $y_{s} <y_{0}$
and $p_{s} \in D$.

First, we show that $\Phi_{t} D\subseteq D$, and hence
$\Phi_{t}\overline{D}\subseteq \overline{D}$ for $t>0$. We
consider the solution $p(t)=(x(t),y(t),z(t))$ with $p(0)\in D$. On
the faces $z=0$ and $z=z_{0}$ of the bar~$D$, we have
$\dot{z}=G(y,0)>0$ and $\dot{z}=-T+G(y,z_{0})<0$, respectively,
and hence $0<z(t)<z_{0}$ for $t>0$. On the faces $x=0$ and
$x=x_{0}$, we have $\dot{x}=R(z)>0$ and $\dot{x}(t)=-1+R(z(t))<0$
for $p(t)$, respectively, and hence $0<x(t)<x_{0}$ for $t>0$. On
the faces $y=0$ and $y=y_{0}$, we respectively have
$\dot{y}(t)=x(t)-G(0,z(t))=x(t)>0$ and
$\dot{y}(t)=x(t)-G(y_{0},z(t))<x_{0}-G(y_{0},0)=0$ for $p(t)$,
whence $0<y(t)<y_{0}$ for $t>0$.

We write $\Pi =\{p\in \partial D:\Phi_{t} p\in D,\,t>0\}$,
$\Pi_{0} =\partial D \setminus \Pi $, and $p_{0} =(x_{0},y_{0}
,0)$. We see that $\Pi_{0} \subseteq l_{1} \bigcup l_{2} \bigcup
l_{3}\bigcup \{ p_{0} \}$, where $l_{1}:\{ x=x_{0},\,0\leq
y<y_{0},\,z=0\}$, $l_{2}:\{ x=0,\, y=0,\,0\leq z \leq z_{0} \}$,
and $l_{3}:\{x=x_{0},\, y=y_{0},\,0\leq z\leq z_{0} \}$. On
$l_{1}$, $l_{2}$, and $l_{3}$, we respectively have $\dot{z}>0$,
$\dot{x}>0$, $\dot{x}<0$, and hence $\Phi_{t} p\in D$, $t>0$, on
$\Pi_{0}\setminus \{p_{0}\}$. Since $\Phi_{t} p_{0} \ne p_{0}$, we
see that $\Phi_{t} p_{0}\in D$, $t>0$. Thus, $\Pi_{0} =\phi$, $\Pi
=\partial D$, and $\Phi_{t} \overline{D}\subseteq D$ for $t>0$.
$\Box$

\noindent \textbf{4.2. Inertial manifold.}
In the natural decomposition $f=-A+F$ of the vector field $f$
of system~(4.1) into the linear and nonlinear parts, we have
\[
A=\left(\begin{array}{ccc} {-k} & {0} & {0} \\ {0} & {0} & {0} \\ {0} & {0} & {-q} \end{array}\right)\;{\kern 1pt},
\quad \;{\kern 1pt}
F\left(\begin{array}{c} {x} \\ {y} \\ {z} \end{array}\right)=\; \,
 \left(\begin{array}{c} {R(z)} \\ {x-G(y,z)} \\ {G(y,z)} \end{array}\right){\kern 1pt} {\kern 1pt} {\kern 1pt} .
\]

This decomposition with symmetric matrix $A$ is chosen in order to
best provide condition (iii) of Theorem 2.1. For the matrix $A$ we
have $\lambda_{1}=0$, $\lambda_{2}=q$, $\lambda_{3}=k$. The change
$u=y+z$ takes (4.1) to the form
$$
\dot{x}=-kx+R(z),\quad \dot{u}=x-qz, \quad \dot{z}=-qz+G(u-z,z)
\eqno (4.2)
$$
in the variables $(x,u,z)$ with the vector field decomposition $f_{1}=-A+F_{1}$,
where $F_{1}:(x,u,z)\rightarrow (R(z),\,x-qz,\,G(u-z,z))$.
In this case,
\[
\left(\begin{array}{c} {x} \\ {y} \\ {z} \end{array}\right)=
C\left(\begin{array}{c} {x} \\ {u} \\ {z} \end{array}\right),
\quad \,
C=\left(\begin{array}{ccc} {1} & {0} & {0} \\ {0} & {1} & {-1} \\ {0} & {0} & {1}
\end{array}\right)\; ,
\quad
C^{-1} =\left(\begin{array}{ccc} {1} & {0} & {0} \\ {0} & {1} & {1} \\ {0} & {0} & {1}
\end{array}\right)\, .
\]
The nonlinear part $F_{1}$ in (4.2) is simpler than the nonlinear
part~$F$ in the original system~(4.1), which allows us to sharpen
the estimate of $K=K(k,q)$ for the norm of its Jacobi matrix in
the spectral gap condition $\lambda_{3} -\lambda_{2} >2K$. The
domain $C^{-1}D$ is strictly positive invariant for (4.2). We put
$$
K=\mathop{\max \;  }\limits_{C^{-1} D} \left\|
F_{1} ^{{'} }(\overline{p})\right\|_{2}=\mathop{\max \;
}\limits_{D} \left\| (F_{1} ^{{'} }
 \, C^{-1} ) (p)\right\|_{2},
\quad
F_{1} ^{{'} } \, C^{-1} =\left(\begin{array}{ccc} {0} & {0} & {-R_{z} } \\ {1} &
{0} & {-q} \\ {0} & {G_{y} } & {G_{z} -G_{y} }
\end{array}\right),
\eqno (4.3)
$$
where $\overline{p}=(x,u,z)$. The condition (2.4) of existence of
the inertial manifold means that (1.2) is satisfied for the
function $F_{1}$ on ${\mathbb R}^{3}$. In this connection, it is
useful to consider a $C^{1+\alpha}$ extension of ~$F_{1}$ from the
domain~$C^{-1} D$ to ${\mathbb R}^{3}$ with the same value of~$K$.
To this end, consider the functions $\overline{R}$ and
$\overline{G}$ defined as follows. The function $\overline{R}$
satisfies $\overline{R}(0)=R(0)$ and its derivative
$\overline{R}_{z}$ is an even $2z_{0}$-periodic extension of
$R_{z}$ from $[0,z_{0} ]$ to $\mathbb{R}$. Similarly,
$\overline{G}$ satisfies $\overline{G}(0,0)=G(0,0)$  and its
derivatives $\overline{G}_{y}$ and $\overline{G}_{z}$ are even,
with respect to both $y$ and $z$, and $(2y_{0},2z_{0})$-periodic
extensions of $G_{y}$ and $G_{z}$, correspondingly, from $[0,y_{0}
]\times [0,z_{0} ]$ to $\mathbb{R}^{2}$. If we now put
$F_{2}:(x,u,z)\rightarrow(\overline{R}(z),\,x-qz,\,\overline{G}(u-z,z))$,
then the function $F_{2}$ yields the sought extension of $F_{1}$
from $C^{-1} D$ to ${\mathbb R}^{3}$. Clearly, the phase dynamics
of system~(4.2) in the domain $C^{-1}D$ remains the same
when~$F_{1}$ is replaced by~$F_{2}$.

Let $\Theta =\{ (k,q)\in {\mathbb R}_{+}^{2} ,\; k-q>2K(k,q)\}$.
Then $\lambda_{3}-\lambda_{2}=k-q$ and, for $(k,q)\in \Theta$, the system of equations
$$
\dot{x}=-kx+\overline{R}(z),\quad \dot{u}=x-qz, \quad \dot{z}=-qz+\overline{G}(u-z,z)
\eqno (4.4)
$$
admits a two-dimensional inertial manifold.
The same is also true for the system
$$
\dot{x}=-kx+\overline{R}(z),\quad \dot{y}=x-\overline{G}(y,z), \quad \dot{z}=-qz+\overline{G}(y,z),
\eqno (4.5)
$$
which inherits the phase dynamics of (4.1) in the domain~$D$.

\noindent \textbf{Remark~4.1.}
If $(k_{0} ,q_{0} )\in \Theta$, then $(k,q)\in \Theta $ for $k\ge k_{0}$, $q\ge q_{0}$,
$k-q\ge k_{0}-q_{0}$.

Indeed, since the strictly positive invariant domain~$D$ decreases
as $k$ and $q$ increase, it follows that the constant $K=K(k,q)$
in~(4.3) does not increase and the inequality $k-q>2K$ still
holds. We see that systems (4.1) and (4.5) demonstrate the
two-dimensional final dynamics in the vast domain $\Theta $ of the
parameters $(k,q)$.

\noindent \textbf{4.3. Instability of the singular point.} The
singular points of systems (4.1) and (4.4) are simultaneously
stable or unstable. The Jacobi matrix $f ^{{'}}(p_{s})$ of the
vector field of system (4.1) at the singular point $p_{s}
=(x_{s},y_{s},z_{s})\in D$ has the form
\[
\left(\begin{array}{ccc} {-k} & {0} & {-b} \\ {1} & {-c} & {-d} \\ {0} & {c} & {d-q} \end{array}\right){\kern 1pt} \,
\]
with $b=-R_{z} (z_{s})$, $c=G_{y} (y_{s},z_{s} )$, and $d=G_{z}(y_{s},z_{s})$.
By the Routh--Hurwitz criterion, this point is asymptotically unstable
if $a_{1}<0$ or $a_{1}a_{2}-a_{3}<0$ or $a_{3}<0$, where
$$
a_{1} =c-d+k+q, \quad a_{2} =k(c-d)+qc+kq,\quad a_{3} =(kq+b)c\,.
$$
Because $a_{3}>0$, the point $p_{s}$ is unstable under the condition $a_{2} <0$.
We have $\mathrm{det}f^{{'}}(p_{s})=c(b-kq)$.

\noindent \textbf{4.4. Stable limit cycle.} The complicated
character of nonlinearity in (4.1) requires the use of
computational tools (Maple package) for estimating the Lipschitz
constant $K(k,q)$ and analyzing the instability of~$p_{s}$. As an
example, we take two pairs of parameters $k>q$ and estimate the
norms for the points $p\in D$. The square numerical matrices~$B$
satisfy the inequality $\left\| B \right\|_{2} \le \sqrt{\left\|
B\right\|_{\infty } \cdot \left\| B\right\|_{1} }$, where
$\left\|B\right\|_{\infty } $ and $\left\| B\right\|_{1}$ are the
norms of the linear operators corresponding to~$B$ in~${\mathbb
R}_{\infty }^{n} $ and ${\mathbb R}_{1}^{n} $.

For $k=3$ and $q=0.1$, we have:
$$
\,y_{1} \approx 186, \; x_{s} \approx 0.117,\; y_{s} \approx
49.653,\; z_{s} \approx 1.167 , \; b-kq \approx 0.480,\;  a_{2}
\approx -0.05,
$$
$$
\left\|(F'_{1}\, C^{-1} )(p)\right\|_{\infty }
\le 1.209, \;\, \left\| (F'_{1} \,\, C^{-1} )(p) \right\|_{1} \le
1.166, \;\, \left\|(F'_{1} \, C^{-1} )(p) \right\|_{2} \le
K=1.187.
$$

For $k=2.5$ and $q=0.1$, we have:
$$
\,y_{1} \approx 204,\; x_{s} \approx 0.123,\; y_{s} \approx
49.558,\; z_{s} \approx 1.230,\; b-kq \approx 0.438,\; a_{2}
\approx -0.01,
$$
$$
\left\|(F'_{1} \, C^{-1} )(p)\right\|_{\infty }
\le 1.209, \; \left\| (F'_{1} \, C^{-1} )(p) \right\|_{1} \le
1.166, \; \left\|(F'_{1}\, C^{-1} )(p)\right\|_{2} \le K=1.187.
$$

The vector field of system (4.4) is real analytic in the strictly
positive invariant domain $C^{-1}D$, and this domain contains a
unique singular point. In both cases $a_{2} <0$, $\mathrm{det}f
^{{'}}(p_{s})=c(b-kq)\neq 0$, and $k-q>2K$, so that by
Theorem~2.1, system~(4.4) admits a stable limit cycle $\Gamma\in
C^{-1}D$ for the chosen values of~$k$ and~$q$. It is easy to trace
the continuous dependence of the quantities $K=K(k,q)$,
$b=b(k,q)$, and $a_{2} =a_{2} (k,q)$ on their arguments, and thus,
the system admits stable periodic regimes for the parameters
$(k,q)$ in sufficiently small neighborhoods of the points
$(3,0.1)$ and $(2.5,0.1)$. This implies that, for the same values
of $(k,q)$, the original system~(4.1) has a stable limit cycle
localized in the domain~$D$.

\bigskip
\noindent \textbf{5. Conclusion}

The spectral gap method is based on the presence of a natural
self-adjoint linear component~$-A$ of the vector field of ODE with
dominating third eigenvalue, $\lambda_{3}(A)>\lambda_{2}(A)$,
which somewhat restricts the range of applications. The advantages
of the method are the transparency of statements and the relative
simplicity of its use. The problems solved by this method are
technically reduced to careful estimation of the Lipschitz
constant in the nonlinear component of the equations and
determination of a strictly positive invariant domain in the phase
space that contains a unique (asymptotically unstable) singular
point. In general, the proposed method can well complement the
list of well-known approaches to the problem of determining stable
limit cycles of ordinary differential equations in~${\mathbb
R}^{n}$, lying on invariant $2\mathrm{D}$-manifolds of the
Cartesian structure.

Existence of an inertial manifold of dimension greater than 2 is
also of interest. For example, the presence of such manifolds of
dimension 3 guarantees, that all invariant tori (if any) of the
dynamical system lie on the invariant three-dimensional
$C^{1}$-manifold of the form (2.1). In the most common spectral
gap condition (2.4) allows us to state that the union of all
bounded invariant sets lies on the smooth invariant
$m$-dimensional manifold of the Cartesian structure.

\bigskip
\noindent \textbf{References}

\noindent [1] G.~Sell and Y.~You, \textit{Dynamics of Evolutionary
Equations}, Springer, New York, 2002.

\noindent [2] S.~Zelik, Inertial manifolds and finite-dimensional
reduction for dissipative PDEs,
\textit{Proc. Roy. Soc. Edinburgh, Ser. A} \textbf{144} (2014), No.~6, 1245--1327.

\noindent [3] M.~Miklavcic, A sharp condition for existence of an
inertial manifolds, \textit{J. Dyn. Differ. Equations} \textbf{3}
(1991), No.~3, 437--456.

\noindent [4] A.~V.~Romanov, Sharp estimates of the dimension of
inertial manifolds for nonlinear parabolic equations, \textit{Izv.
Math.} \textbf{43} (1994), No.~1, 31--47.

\noindent [5] J.~A.~Yorke, Periods of periodic solutions and
Lipschitz constant, \textit{Proc. Amer. Math. Soc.} \textbf{22}
(1969), No.~2, 509--512.

\noindent [6] M.~A.~Krasnosel'skii and P.~P.~Zabreiko,
\textit{Geometrical Methods of Nonlinear Analysis}, 2nd ed.,
Springer, New York, 2011.

\noindent [7] R.~A.~Smith, Certain differential equations have
only isolated periodic orbits, \textit{Ann. Mat. Pura Appl.},
\textbf{137} (1984), No.~1, 217--244.

\noindent [8] R.~Reissig, G.~Sansone and R.~Conti,
\textit{Qualitative Theorie Nichtlinearer}

\noindent \textit{Differentialgleichungen}, Edizioni Cremonese,
Roma, 1963.

\noindent [9] R.~A.~Smith, Poincare index theorem concerning
periodic orbits of differential equations, \textit{Proc. London
Math. Soc.} (3) \textbf{48} (1984), No.~2, 341--362.

\noindent [10] R.~A.~Smith, Orbital stability for ordinary
differential equations, \textit{J. Differ. Equations} \textbf{69}
(1987), No.~2, 265--287.

\noindent [11] V. V. Sidorenko and A. Celletti, “Spring–mass”
model of tethered satellite systems: properties of planar periodic
motions, \textit{Celestial Mechanics and Dynamical Astronomy},
\textbf{107} (2010), No.~1--2, 209–-231.

\noindent [12] I.~A.~Galiullin and L.~A.~Kondratieva, Satellite
inertial manifolds and boundary cycles, \textit{Cosmonautics and
Rocket Engineering} (2011), No.~3, 73--76 (in Russian).

\noindent [13] C.~Suguna, K.~K.~Chowdhury and S.~Sinha, Minimal
model for complex dynamics in cellular processes, \textit{Phys.
Rev.~E} \textbf{60} (1999), No.~5, 5943--5949.

\end{document}